\documentclass{article}
\usepackage{amssymb}

\usepackage{amsmath}

\author{Gert Almkvist}
\title{Asymptotics of various partitions}

\begin{document}

\maketitle

\textbf{1. Introduction.}

In [1] the author gave the first two terms in the asymptotic expansion of $%
\widetilde{p}(r,n)=$ the number of partitions of $n$ into parts $\geq r$.
In this paper we compute further terms until we get all digits of $%
\widetilde{p}(12,1200).$

Further examples treated are ''basic partitions'' ( Nolan, Savage, Wilf [7]),
3-colored partitions ( Kolitsch [5]), plane partitions strictly decreasing
along rows ( Macdonald [6] ), the number of semisimple p-rings (Knopfmacher
[4] ) and concave partitions (Snellman,Paulsen [8]). 
\[
\]

\textbf{2. Partitions without small parts.}

We have the generating function 
\[
\widetilde{f}_r(x)=\sum_{n=0}^\infty \widetilde{p}(r,n)x^n=\frac{F(x)}{%
f_{r-1}(x)} 
\]
where 
\[
F(x)=\sum_{n=0}^\infty p(n)x^n=\prod_{j=1}^\infty (1-x^j)^{-1} 
\]
and 
\[
f_r(x)=\prod_{j=1}^r(1-x^j)^{-1} 
\]
\[
\]
In [1] we used derivatives of 
\[
\sinh (\sqrt{a\xi }) 
\]
to express the formulas. Here we will rather use associated Bessel functions
with half integer index.By Example 5 in [1] we have 
\[
\widetilde{f}_r(e^{-t})=\frac{(r-1)!}{\sqrt{2\pi }}t^{r-1/2}\exp \left\{
\sum_{j=1}^\infty \frac{B_{2j}B_{2j+1}(r)}{2j(2j+1)!}t^{2j}\right\} \exp
\left\{ \varsigma (2)t^{-1}-\frac{B_2(r)t}4\right\} 
\]
\[
\]
Remember that formally we have 
\[
t^{-a}\rightarrow \frac{n^a}{\Gamma (\alpha )} 
\]
Hence 
\[
t^a\exp (ct^{-1}+bt)=e^{bt}\sum_{j=0}^\infty \frac{c^j}{j!}%
t^{-j+a}\rightarrow 
\]
\[
\sum_{j=0}^\infty \frac{c^j\xi ^{j-a-1}}{j!\Gamma (j-a)}=(\frac c\xi
)^{(1+a)/2}I_{-1-a}(2\sqrt{c\xi }) 
\]
where 
\[
\xi =n+b 
\]
Hence if 
\[
\exp \left\{ \sum_{j=1}^\infty \frac{B_{2j}B_{2j+1}(r)}{2j(2j+1)!}%
t^{2j}\right\} =\sum_{j=0}^\infty c_{2j}t^{2j} 
\]
then 
\[
\Phi _1(n)=\frac{(r-1)!}{\sqrt{2\pi }}\sum_{j=0}^\infty c_{2j}\left( \frac{%
\varsigma (2)}\xi \right) ^{-r/2+1/4+j}I_{-(r+1/2+2j)}(2\sqrt{\varsigma
(2)\xi }) 
\]
where 
\[
\xi =n-\frac{B_2(r)}4 
\]
\[
\]

For the second approximation belonging to the singular point $x=-1$ of $%
\widetilde{f}_r(x)$ we obtain for even $r$%
\[
\widetilde{f}_r(-e^{-t})=\frac{2^{r-1}(\frac r2-1)!}{\sqrt{\pi }}%
t^{(r-1)/2}\exp \left\{ 
\begin{array}{c}
\frac{\varsigma (2)t^{-1}}4-\frac{B_2(r)t}4+\sum_{j=1}^\infty \frac{%
2^{2j}B_{2j}t^{2j}}{2j(2j+1)!}(B_{2j+1}(\frac r2) \\ 
+(2^{2j}-1)B_{2j+1}(\frac{r+1}2)
\end{array}
\right\} 
\]
\[
\]
Assuming that the last part of the exponential function is 
\[
\sum_{j=0}^\infty d_{2j}t^{2j} 
\]
we find 
\[
\Phi _2(n)=(-1)^n\frac{2^{r-1}(\frac r2-1)!}{\sqrt{\pi }}\sum_{j=0}^\infty
d_{2j}\left( \frac{\varsigma (2)}{4\xi }\right)
^{r/4+1/4+j}I_{-(r/2+1/2+2j)}(\sqrt{\varsigma (2)\xi }) 
\]
where 
\[
\xi =n-\frac{B_2(r)}4 
\]
\begin{eqnarray*}
&& \\
&&
\end{eqnarray*}

\textsl{Higher order terms}.

We need to approximate $\widetilde{f}_r(x)$ near $x=\exp (2\pi ih/k)$ where $%
(h,k)=1.$ This will lead to some formidable expressions containing Dedekind
sums and various other sums of derivatives of $\cot (x)$ evaluated at $%
x=jh\pi /k.$ We will save the reader from this. Instead we show how to do it
numerically in an example. We continue Example 5 in [1]. Thus let $r=12$ and 
$n=1200.$ Take $k=7$ which is the first really difficult case. We consider
the following three singular points:

1. $x=\exp (2\pi i/7).$

In general we have 
\[
F(\exp (2\pi ih/k-t))\thickapprox \sqrt{\frac{kt}{2\pi }}\exp \left\{ \frac{%
\varsigma (2)t^{-1}}{k^2}-\frac t{24}+\pi is(h,k)\right\} 
\]
where $(h,k)=1.$ Here 
\[
s(h,k)=\sum_{j=1}^{k-1}((jh/k))((j/k)) 
\]
where 
\[
((x))=x-[x]-\frac 12\text{ if }x\notin \mathbf{Z} 
\]
\[
=0\text{ otherwise} 
\]
Then 
\[
F(\exp (2\pi i/7-t))\thickapprox \sqrt{\frac{7t}{2\pi }}\exp \left\{ \frac{%
\varsigma (2)}{49t}-\frac t{24}+\frac{5\pi i}{14}\right\} 
\]
since 
\[
s(1,7)=\frac 5{14} 
\]
\textsl{Warning. }

In the beginning I used the table in Hardy-Ramanujan's famous paper [3] ,
but there are a lot of errors (even in the Chelsea edition of Ramanujan's
collected papers ). So it is better to do your own computations. 
\[
\]

2. $x=\exp (4\pi i/7).$

Then 
\[
F(\exp (4\pi i/7-t))\thickapprox \sqrt{\frac{7t}{2\pi }}\exp \left\{ \frac{%
\varsigma (2)}{49t}-\frac t{24}+\frac{\pi i}{14}\right\} 
\]
since 
\[
s(2,7)=\frac 1{14} 
\]

3. $x=\exp (6\pi i/7)$

Then 
\[
F(\exp (2\pi i/7-t))\thickapprox \sqrt{\frac{7t}{2\pi }}\exp \left\{ \frac{%
\varsigma (2)}{49t}-\frac t{24}-\frac{\pi i}{14}\right\} 
\]
since 
\[
s(3,7)=-\frac 1{14} 
\]

The other points with $h=4,5,6$ are treated by conjugation. We have 
\[
\widetilde{f}_r(x)=g(x)F(x) 
\]
where 
\[
g(x)=\prod_{j=1}^{11}(1-x^j)=1-x-x^2+x^5+x^7-2x^{15}+\cdot \cdot \cdot
x^{65}-x^{66} 
\]
We expand 
\[
g(\exp (2\pi i/7-t))=\sum_{j=0}^\infty c_jt^j 
\]
numerically. Thus 
\[
c_0=0 
\]
\[
c_1=-56.249-246.445i 
\]
\[
c_2=3157.20+7835.75i 
\]
\[
etc 
\]
Put 
\[
\xi =1200-\frac 1{24} 
\]
Then the contribution from $x=\exp (2\pi i/7)$ and $x=\exp (-2\pi i/7)$ is 
\[
\Phi _{7,1}=2\sqrt{\frac 7{2\pi }}\sum_{j=1}^{10}Re\left\{ c_j\exp \left\{ 
\frac{5\pi i}{14}-\frac{2\pi i1200}7\right\} \right\} \left\{ \frac{%
\varsigma (2)}{49\xi }\right\} ^{3/4+j/2}I_{-3/2-j}(2\sqrt{\frac{\varsigma
(2)\xi }{49}}) 
\]
\[
=-26.49890 
\]
Similarly we obtain 
\[
\Phi _{7,2}=9.34180 
\]
\[
\Phi _{7,3}=8.86403 
\]
and 
\[
\Phi _7=\Phi _{7,1}+\Phi _{7,2}+\Phi _{7,3}=-8.29306 
\]
Taking 16 terms we get 
\[
\Phi _1=4\text{ }90015\text{ }90791\text{ }72948\text{ }38443\text{ }%
67842.13295 
\]
\[
\Phi _2=33\text{ }28709\text{ }71879.70797 
\]
\[
\Phi _3=124\text{ }91942.03287 
\]
\[
\Phi _4=51637.99481 
\]
\[
\Phi _5=661.16304 
\]
\[
\Phi _6=164.15289 
\]
\[
\Phi _7=-8.29306 
\]
\[
\Phi _8=3.22670 
\]
\[
\Phi _9=-0.19405 
\]
\[
\Phi _{10}=0.23922 
\]
\[
\Phi _{11}=0.00072 
\]
\[
\Phi _{12}=0.34102 
\]
\[
\Phi _{13}=0.89399 
\]
\[
\Phi _{14}=-0.64967 
\]
\[
\Phi _{15}=0.94330 
\]
\[
\Phi _{16}=0.44736 
\]
Summing up 
\[
\sum_{j=1}^{16}\Phi _j=4\text{ }90015\text{ }90791\text{ }72981\text{ }67278%
\text{ }84124.17000 
\]
which agrees with 
\[
\widetilde{p}(12,1200)=4\text{ }90015\text{ }90791\text{ }72981\text{ }67278%
\text{ }84124 
\]
\[
\]
\textbf{Remark.}

We have $\Phi _{11}=0.00072$ but from $k=12$ the $\Phi _k$ become much
larger. This depends on the fact that the constant term 
\[
\prod_{j=1}^{11}(1-\exp (2\pi ijh/k)) 
\]
is no longer zero. 
\[
\]

\textbf{3. Basic partitions.}

Our next example concerns a limiting case of''basic partitions'' (see
Nolan,Savage,Wilf [7] ). Thus we have the generating function 
\[
f(x)=\prod_{j=1}^\infty \frac{1+x^j}{1-x^j}=\sum_{n=0}^\infty a_nx^n 
\]
Then 
\[
f(x)=\prod_{j=1}^\infty \frac{1-x^{2j}}{(1-x^j)^2}=\frac{F(x)^2}{F(x^2)} 
\]
and 
\[
f(\exp (-t))=\frac{F(\exp (-t))^2}{F(\exp (-2t))}\thickapprox \frac 1{2\sqrt{%
\pi }}t^{1/2}\exp (\frac{\pi ^2}{4t}) 
\]
which gives 
\[
\Phi _1(n)=\frac 1{2\sqrt{\pi }}\left( \frac an\right) ^{3/4}I_{-3/2}(2\sqrt{%
an}) 
\]
where 
\[
a=\frac{\pi ^2}4 
\]
\textsl{General k.}

Only $x=\exp (2\pi ih/k)$ with odd $k$ are singular points of $f.$ We have 
\[
f(\exp (2\pi ih/k-t))=\frac{F(\exp (2\pi ih/k-t))^2}{F(\exp (4\pi ih/k-2t))} 
\]
\[
\thickapprox \frac{\sqrt{kt}}{2\sqrt{\pi }}\exp \left\{ \pi
i(2s(h,k)-s(2h,k))\right\} \exp \left\{ \frac a{k^2t}\right\} 
\]
Put 
\[
T(h,k)=2s(h,k)-s(2h,k) 
\]
and 
\[
\widetilde{A}(k,n)=\sum_{(h,k)=1}\exp \left\{ \pi iT(h,k)-2\pi inh/k\right\} 
\]
Then 
\[
\Phi _k(n)=\frac 12\sqrt{\frac k\pi }\left( \frac a{k^2n}\right) ^{3/4}%
\widetilde{A}(k,n)I_{-3/2}(2\sqrt{\frac{an}{k^2}}) 
\]
where

$\widetilde{A}(1,n)=1$

$\widetilde{A}(3,n)=2\cos (2\pi n/3-\pi /6)$

$\widetilde{A}(5,n)=2\left\{ \cos (2\pi n/5-2\pi /5)+\cos (4\pi n/5-\pi
/5)\right\} $

$\widetilde{A}(7,n)=2\left\{ \cos (2\pi n/7-9\pi /14)+\cos (4\pi n/7-\pi
/14)+\cos (6\pi n/7-3\pi /14)\right\} $

$\widetilde{A}(9,n)=2\left\{ \cos (2\pi n/9-8\pi /9)+\cos (4\pi n/9-4\pi
/9)+\cos (8\pi n/9-2\pi /9)\right\} $

$\widetilde{A}(11,n)=2\left\{ 
\begin{array}{c}
\cos (2\pi n/11-25\pi /22)+\cos (4\pi n/11-7\pi /22)+\cos (6\pi n/11-\pi /22)
\\ 
+\cos (8\pi n/11-9\pi /22)+\cos (10\pi /11-5\pi /22)
\end{array}
\right\} $

$\widetilde{A}(13,n)=2\left\{ 
\begin{array}{c}
\cos (2\pi n/13-18\pi /13)+\cos (4\pi n/13-9\pi /13)+ \\ 
\cos (6\pi n/13-6\pi /13)+\cos (8\pi n/13+2\pi /13)+ \\ 
\cos (10\pi n/13-\pi /13)+\cos (12\pi n/13-3\pi /13)
\end{array}
\right\} $%
\[
\]
Example. Let $n=1000.$ Then 
\[
a_{1000}=17293\text{ }58313\text{ }74933\text{ }37582\text{ }44155\text{ }%
69812\text{ }30246\text{ }17584 
\]
\[
\]
We have

$\Phi _1=17293$ $58213$ $74933$ $37582$ $44155$ $69812$ $30245$ $67360.51591$

$\Phi _3=0$

$\Phi _5=50530.87514$

$\Phi _7=-305.16980$

$\Phi _9=0$

$\Phi _{11}=-1.78371$

$\Phi _{13}=-0.43964$

\[
\]
\[
\Phi _1+\cdot \cdot \cdot \Phi _{13}=17293\text{ }58213\text{ }74933\text{ }%
37582\text{ }44155\text{ }69812\text{ }30246\text{ }17583.99790 
\]
so the error is $-0.0021.$%
\[
\]

\textbf{4. 3-colored partitions.}

We have the generating function 
\[
f(x)=\prod_{j=1}^\infty \frac{(1-x^{9j})^3}{(1-x^{3j})(1-x^j)^3}%
=\sum_{n=0}^\infty a_nx^n 
\]
This case caused me a lot of trouble until I realized that one has to use
the \textsl{exact }functional equation for 
\[
F(x)=\prod_{j=1}^\infty (1-x^j)^{-1} 
\]
namely 
\[
F(\exp (-t))=\sqrt{\frac t\pi }\exp \left\{ \frac{\pi ^2}{6t}-\frac
t{24}\right\} F(\exp \left\{ -\frac{4\pi ^2}t\right\} 
\]
Now 
\[
f(x)=\frac{F(x)^3F(x^3)}{F(x^9)^3} 
\]
gives 
\[
f(\exp (-t))=\sqrt{\frac t{486\pi }}\exp \left\{ \frac{\pi ^2}{2t}+\frac{7t}%
8\right\} \left\{ 1-3\exp \left\{ -\frac{4\pi ^2}{9t}\right\} +\cdot \cdot
\cdot \right\} 
\]
\[
\thickapprox \frac 1{\sqrt{486\pi }}\left\{ t^{1/2}\exp \left\{ \frac{\pi ^2%
}{2t}+\frac{7t}8\right\} -3t^{1/2}\exp \left\{ \frac{\pi ^2}{18t}+\frac{7t}%
8\right\} \right\} 
\]
\[
\]
The other terms coming from $F(\exp (-4\pi ^2/t))$ and $F(\exp (-4\pi
^2/3t)) $ give terms of the form $\exp (-c/t)$ where $c>0$ which transform
into ordinary Bessel functions $J_{-3/2}(2\sqrt{c\xi })$ which are very
small for large $\xi .$

We get as usual with $\xi =n+7/8$%
\[
\Phi _1(n)=\frac 1{\sqrt{486\pi }}\left\{ \left( \frac{\pi ^2}{2\xi }\right)
^{3/4}I_{-3/2}(2\sqrt{\frac{\pi ^2\xi }2})-3\left( \frac{\pi ^2}{18\xi }%
\right) ^{3/4}I_{-3/2}(2\sqrt{\frac{\pi ^2\xi }2)}\right\} 
\]
More generally we have 
\[
F(\exp (2\pi ih/k-t))= 
\]
\[
\sqrt{\frac{kt}{2\pi }}\exp \left\{ \pi is(h,k)+\frac{\pi ^2}{6k^2t}-\frac
t{24}\right\} F(\exp (\left\{ 2\pi ih^{\prime }/k-\frac{4\pi ^2}{k^2t}%
\right\} ) 
\]
where 
\[
hh^{\prime }\equiv -1\text{ (}mod\text{ }k) 
\]
In particular for $k=2$ we get 
\[
f(-\exp (-t))\thickapprox \sqrt{\frac t{243\pi }}\left\{ \exp \left\{ \frac{%
\pi ^2}{8t}+\frac{7t}8\right\} +3\exp \left\{ \frac{\pi ^2}{72t}+\frac{7t}%
8\right\} \right\} 
\]
and thus 
\[
\Phi _2(n)=(-1)^n\frac 1{\sqrt{243\pi }}\left\{ \left( \frac{\pi ^2}{8\xi }%
\right) ^{3/4}I_{-3/2}(2\sqrt{\frac{\pi ^2\xi }8})+3\left( \frac{\pi ^2}{%
72\xi }\right) ^{3/4}I_{-3/2}(2\sqrt{\frac{\pi ^2\xi }{72}})\right\} 
\]
Further computations give (one has to be careful when $k$ is divisible by $3$
).

\[
\Phi _3(n)=2\cos (\frac{2\pi n}3-\frac \pi 6)\frac 1{\sqrt{18\pi }}\left( 
\frac{\pi ^2}{18\xi }\right) ^{3/4}I_{-3/2}(2\sqrt{\frac{\pi ^2\xi }{18}}) 
\]
\[
\Phi _4(n)=2\cos (\frac{2\pi n}4+\frac \pi 8)\sqrt{\frac 4{484\pi }}\left( 
\frac{\pi ^2}{32\xi }\right) ^{3/4}I_{-3/2}(2\sqrt{\frac{\pi ^2\xi }{32}}) 
\]
\[
\Phi _5(n)=2\left\{ \cos (\frac{2\pi n}5-\frac{6\pi }5)+\cos (\frac{4\pi n}%
5-\frac \pi 5)\right\} \sqrt{\frac 5{486\pi }}\left( \frac{\pi ^2}{50\xi }%
\right) ^{3/4}I_{-3/2}(2\sqrt{\frac{\pi ^2\xi }{50}}) 
\]
\[
\Phi _6(n)=2\cos (\frac{2\pi n}6-\frac{5\pi }6)\frac 1{\sqrt{9\pi }}\left( 
\frac{\pi ^2}{72\xi }\right) ^{3/4}I_{-3/2}(2\sqrt{\frac{\pi ^2\xi }{72}}) 
\]
\[
\Phi _7(n)=2\left\{ \cos (\frac{2\pi n}7-\frac{11\pi }{14})+\cos (\frac{4\pi
n}7-\frac \pi {14})+\cos (\frac{6\pi n}7-\frac{13\pi }{14})\right\} 
\]
\[
\sqrt{\frac 7{486\pi }}\left( \frac{\pi ^2}{98\xi }\right) ^{3/4}I_{-3/2}(2%
\sqrt{\frac{\pi ^2\xi }{98}}) 
\]
Numerical example.

Take $n=200.$ Then 
\[
a_{200}=1747\text{ }47949\text{ }05123\text{ }77771\text{ }22300 
\]
and

$\Phi _1=1747$ $47949$ $05123$ $51995$ $85303.8855$

$\Phi _2=25778$ $34263.0509$

$\Phi _3=-2$ $97748.8532$

$\Phi _4=474.3640$

$\Phi _5=0$

$\Phi _6=5.5359$

$\Phi _7=0.5261$

Thus 
\[
\Phi _1+\cdot \cdot \cdot \Phi _7=1747\text{ }47949\text{ }05123\text{ }77771%
\text{ }22298.5093 
\]
so the error is $-1.4907$%
\[
\]

\textbf{5.Plane partitions strictly decreasing along rows.}

We have the generating function (see Macdonald [6] p.83) 
\[
f(x)=\prod_{j=1}^\infty (1-x^j)^{-[(j+1)/2]}=\sum_{n=0}^\infty a_nx^n 
\]
Consider 
\[
g(t)=\log (f(\exp (-t)))=-\sum_{j=1}^\infty j\left\{ \log (1-\exp
(-(2j-1)t))+\log (1-\exp (-2jt))\right\} 
\]
\[
=\sum_{j=1}^\infty \sum_{k=1}^\infty \frac jk\left\{ \exp (-k(2j-1)t)+\exp
(-2jkt)\right\} 
\]
We form the Mellin transformation of $g(t)$%
\[
\widetilde{g}(s)=\int_0^\infty g(t)t^{-s}dt= 
\]
\[
\Gamma (s)\sum_{j=1}^\infty \sum_{k=1}^\infty \left\{ \frac
j{k^{1+s}(2j-1)^s}+\frac j{k^{1+s}(2j)^s}\right\} = 
\]
\[
\frac 12\Gamma (s)\varsigma (1+s)\left\{ \varsigma (s-1)+(1-2^{-s})\varsigma
(s)\right\} 
\]
To compute the inverse Mellin transformation 
\[
g(t)=\frac 1{2\pi i}\int_{c-i\infty }^{c+i\infty }\widetilde{g}(s)t^{-s}ds 
\]
$(c>2)$ we need the residues

\[
Res(\widetilde{g}(s)t^{-s},s=2)=\frac{\varsigma (3)}{2t^2} 
\]
\[
Res(\widetilde{g}(s)t^{-s},s=1)=\frac{\pi ^2}{24t} 
\]
\[
Res(\widetilde{g}(s)t^{-s},s=0)=\frac 12\left\{ \varsigma ^{\prime }(-1)-%
\frac{\log (2)}2+\frac{\log (t)}{12}\right\} 
\]
\[
Res(\widetilde{g}(s)t^{-s},s=-1)=\frac t{48} 
\]
\[
Res(\widetilde{g}(s)t^{-s},s=1-2k)=0\text{ for }k\geq 1 
\]
\[
Res(\widetilde{g}(s)t^{-s},s=-2k)=\frac 12\varsigma (1-2k)\varsigma (-1-2k)%
\frac{t^{2k}}{(2k)!}\text{ for }k\geq 1 
\]
Hence 
\begin{eqnarray*}
&&f(\exp (-t)) \\
&=&2^{-1/4}\exp (\varsigma ^{\prime }(-1)/2)t^{1/24}\exp \left\{ \frac{%
\varsigma (3)}{2t^2}+\frac{\pi ^2}{24t}+\frac t{48}+\frac
12\sum_{k=1}^\infty \frac{\varsigma (1-2k)\varsigma (-1-2k)t^{2k}}{(2k)!}%
\right\}
\end{eqnarray*}

Use the notation 
\[
a=\frac{\varsigma (3)}2 
\]
\[
b=\frac{\pi ^2}{24} 
\]
\[
c=2^{-1/4}\exp (\varsigma ^{\prime }(-1)/2) 
\]
\[
\]
\textbf{Remark.}

$\varsigma ^{\prime }(-1)$ is known as \textsl{Kinkelins constant }(see [2]
for some history )

First we replace 
\[
\exp \left\{ \frac 12\sum_{k=1}^\infty \frac{\varsigma (1-2k)\varsigma
(-1-2k)t^{2k}}{(2k)!}\right\} =1+\sum_{k=1}^\infty c_kt^k 
\]
by $1.$ We also disregard $\exp (t/48)$ which is merely the translation $%
n\rightarrow n+1/48=\xi .$%
\[
f(\exp (-t))_0\thickapprox ct^{1/24}\exp (at^{-2}+bt^{-1})= 
\]
\[
c\sum_{j=0}^\infty \sum_{k=0}^\infty \frac{a^jb^kt^{-2j-k+1/24}}{j!k!}%
\rightarrow c\sum_{j=0}^\infty \sum_{k=0}^\infty \frac{a^jb^k\xi
^{2j+k-1/24-1}}{j!k!\Gamma (2j+k-1/24)} 
\]
Hence 
\[
\Phi _{1,0}(n)=c\left( \frac b\xi \right) ^{25/48}\sum_{j=0}^\infty \frac
1{j!}\left( \frac{a\xi }b\right) ^jI_{2j-25/24}(2\sqrt{b\xi }) 
\]
We also have to consider terms coming from $c_2t^2+c_4t^4+\cdot \cdot \cdot $
. We have 
\[
c_2=-\frac 1{5760} 
\]
\[
c_4=-\frac{313}{464486400} 
\]
\[
c_6=-\frac{91207}{8026324992000} 
\]
\[
etc 
\]
The term coming from $c_2t^2$ will be 
\[
\Phi _{1,2}(n)=-\frac 1{5760}c_2\left( \frac b\xi \right)
^{25/48+1}\sum_{j=0}^\infty \frac 1{j!}\left( \frac{a\xi }b\right)
^jI_{2j-25/24-2}(2\sqrt{b\xi }) 
\]
Let us take $n=200.$ Then 
\[
a_{200}=233\text{ }18651\text{ }62179\text{ }63536\text{ }57014 
\]
Taking eight terms we get 
\[
\Phi _1=\Phi _{1,0}+\Phi _{1,2}+\cdot \cdot \cdot \Phi _{1,14}= 
\]
\[
233\text{ }18651\text{ }62179\text{ }61259\text{ }31636.4295 
\]
so we get 14 correct digits out of 23. There is a lot of cancellation so we
have to compute with many more digits than 23. 
\[
\]

\textsl{Second term.}

Let 
\[
g_2(t)=\log (f(-\exp (-t))= 
\]
\[
\sum_{j=1}^\infty j\left\{ \log (1-\exp (-(2j-1)t))-\log (1-\exp
(-2(2j-1)t))-\log (1-\exp (-2jt))\right\} 
\]
with Mellin transformation 
\[
\widetilde{g}_2(s)=\frac 12\Gamma (s)\varsigma (1+s)\left\{ (5\cdot
2^{-s}-2\cdot 2^{-2s}-1)\varsigma (s-1)-(1-2^{-s})\varsigma (s)\right\} 
\]
with residues 
\[
Res(\widetilde{g}_2(s)t^{-s},s=2)=\frac{\varsigma (3)}{16t^2} 
\]
\[
Res(\widetilde{g}_2(s)t^{-s},s=1)=-\frac{\pi ^2}{48t} 
\]
\[
Res(\widetilde{g}_2(s)t^{-s},s=0)=\varsigma ^{\prime }(-1)+\frac{\log (2)}{24%
}+\frac{\log (t)}{12} 
\]
\[
Res(\widetilde{g}_2(s)t^{-s},s=-1)=\frac t{48} 
\]
\[
Res(\widetilde{g}_2(s)t^{-s},s=-2)=-\frac{13t^2}{5760} 
\]
\[
Res(\widetilde{g}_2(s)t^{-s},s=-4)=\frac{433t^4}{1451520} 
\]
\[
Res(\widetilde{g}_2(s)t^{-s},s=-6)=\frac{7873t^6}{97091200} 
\]
\[
etc 
\]
Hence 
\[
f(-\exp (-t))=2^{1/24}\exp (\varsigma ^{\prime }(-1))t^{1/12}\exp \left\{ 
\frac{\varsigma (3)}{16t^2}-\frac{\pi ^2}{48t}+\frac t{48}\right\} 
\]
\[
\exp \left\{ -\frac{13t^2}{5760}+\frac{433t^4}{1451520}+\frac{7873t^6}{%
87091200}+\cdot \cdot \cdot \right\} 
\]
We get 
\[
\Phi _{2,0}(n)=(-1)^n2^{1/24}\exp (\varsigma ^{\prime }(-1))\left( \frac{b_2}%
\xi \right) ^{13/24}\sum_{j=0}^\infty \frac 1{j!}\left( \frac{a_2\xi }{b_2}%
\right) ^jI_{-2j-13/12}(2\sqrt{b_2\xi } 
\]
where 
\[
a_2=\frac{\varsigma (3)}{16} 
\]
\[
b_2=-\frac{\pi ^2}{48} 
\]
\[
\xi =n+\frac 1{48} 
\]
We obtain $\Phi _{2,k}$ for $k=2,4,6,\cdot \cdot \cdot $ by expanding 
\[
\exp \left\{ -\frac{13t^2}{5760}+\cdot \cdot \cdot \right\} 
\]
The result is 
\[
\Phi _2=2281\text{ }49810.1011 
\]
and 
\[
\Phi _1+\Phi _2=233\text{ }18651\text{ }62179\text{ }63540\text{ }81446.5432 
\]
so we get 17 correct digits out of 23. 
\[
\]

\textbf{6. Semisimple p-rings.}

Let $s(n)=$ the number of isomorphy classes of semisimple p-rings of order $%
p^n$ ($p=$prime, see Knopfmacher [4], p.64). We have the generating function 
\[
f(x)=\sum_{n=0}^\infty s(n)x^n=\prod_{j=1}^\infty \prod_{k=1}^\infty
(1-x^{jk^2})^{-1} 
\]
To get the first approximation we consider 
\[
g(t)=\log (f(\exp (-t)))=\sum_{j=1}^\infty \sum_{k=1}^\infty
\sum_{r=1}^\infty \frac{\exp (-jk^2rt)}r 
\]
The Mellin transformation is 
\[
\widetilde{g}(s)=\int_0^\infty g(t)t^{s-1}dt=\Gamma (s)\varsigma
(1+s)\varsigma (s)\varsigma (2s) 
\]
with residues 
\[
Res(\widetilde{g}(s)t^{-s},s=1)=\frac{\varsigma (2)^2}t 
\]
\[
Res(\widetilde{g}(s)t^{-s},s=1/2)=\frac 12\sqrt{\frac \pi t}\varsigma
(1/2)\varsigma (3/2) 
\]
\[
Res(\widetilde{g}(s)t^{-s},s=0)=\frac{3\log (2\pi )}4-\frac{\log (t)}4 
\]
Hence 
\[
g(t)=\frac{3\log (2\pi )}4-\frac{\log (t)}4+\frac{\varsigma (2)^2}t+\frac{%
\sqrt{\pi }}2\varsigma (1/2)\varsigma (3/2)t^{-1/2} 
\]
Let 
\[
a=\varsigma (2)^2 
\]
\[
b=\frac{\sqrt{\pi }}2\varsigma (1/2)\varsigma (3/2) 
\]
Then 
\[
f(\exp (-t))=(2\pi )^{3/4}t^{-1/4}\exp (at^{-1}+bt^{-1/2})= 
\]
\[
(2\pi )^{3/4}\sum_{j=0}^\infty \sum_{k=0}^\infty \frac{a^jb^kt^{-j-k/2-1/4}}{%
j!k!}\rightarrow 
\]
\[
(2\pi )^{3/4}\sum_{j=0}^\infty \sum_{k=0}^\infty \frac{a^jb^kn^{j+k/2+1/4-1}%
}{j!k!\Gamma (j+k/2+1/4)} 
\]
Hence 
\[
\Phi _1(n)=(2\pi )^{3/4}\left( \frac an\right) ^{3/8}\sum_{k=0}^\infty \frac
1{k!}\left( b\left( \frac na\right) ^{1/4}\right) ^kI_{k/2-3/4}(2\sqrt{an}) 
\]
To get the second term $\Phi _2(n)$ we consider 
\[
g_2(t)=\log (f(\exp (-t)))=-\sum_{j=0}^\infty \sum_{k=0}^\infty \log
(1-(-1)^{jk^2}\exp (-jk^2t)) 
\]
After some computations one obtains 
\[
\widetilde{g}_2(s)=\Gamma (s)\varsigma (1+s)\varsigma (s)\left\{ -1+3\cdot
2^{-s}+2^{-2s}-3\cdot 2^{-3s}+2^{-4s}\right\} 
\]
with residues 
\[
Res(\widetilde{g}_2(s)t^{-s},s=1)=\frac{7\varsigma (2)^2}{16t}=\frac{a_2}t 
\]
\[
Res(\widetilde{g}_2(s)t^{-s},s=1/2)=\frac 18\sqrt{\frac \pi {8t}}\left\{ 3%
\sqrt{2}-1\right\} \varsigma (1/2)\varsigma (3/2)=\frac{b_2}{\sqrt{t}} 
\]
\[
Res(\widetilde{g}_2(s)t^{-s},s=0)=\frac{3\log (2\pi )}4-\frac{\log (t)}4 
\]
so 
\[
f(-\exp (-t))=(2\pi )^{3/4}t^{-1/4}\exp \left\{
a_2t^{-1}+b_2t^{-1/2}\right\} 
\]
which is exactly $f(\exp (-t))$ with $a\rightarrow a_2$ and $b\rightarrow
b_2.$ Hence 
\[
\Phi _2(n)=(-1)^n(2\pi )^{3/4}\left( \frac{a_2}n\right)
^{3/8}\sum_{j=0}^\infty \left( b_2\left( \frac n{a_2}\right) ^{1/4}\right)
^kI_{k/2-3/4}(2\sqrt{a_2n}) 
\]

\textsl{Numerical Example.}

Let $n=200.$ Then 
\[
s(200)=26122\text{ }95856\text{ }86401 
\]
while

$\Phi _1(200)=26122$ $95856$ $68813.6838$

$\Phi _2(200)=290$ $18442.4542$%
\[
\Phi _1+\Phi _2=26122\text{ }95852\text{ }86711.2679 
\]
so we get 9 correct digits out of 15. The largest terms in $\Phi _1$ are of
size $10^{22}$ so there is a lot of cancellation. Knopfmacher gives the
estimate 
\[
\log (s(n))\thicksim \frac{\pi ^2}3\sqrt{n} 
\]
We can improve this by using the asymptotic formula for Bessel functions. We
have for large $z$%
\[
I_r(z)\thicksim \frac{\exp (z)}{\sqrt{2\pi z}}\left\{ 1-\frac{4r^2-1}{8z}+%
\frac{(4r^2-1)(4r^2-9)}{2!(8z)^2}-\cdot \cdot \cdot \right\} 
\]
Putting 
\[
z=2\sqrt{an} 
\]
$\infty $%
\[
r=\frac k2-\frac 34 
\]
and summing over $k$ we obtain 
\[
s(n)\thicksim \frac{(2\pi \sqrt{a})^{1/4}}{\sqrt{2}n^{5/8}}\exp \left\{ 2%
\sqrt{an}+b\left( \frac na\right) ^{1/4}-\frac{b^2}{16a}\right\} 
\]
For $n=200$ the right hand side is $0.28\cdot 10^{15}$ so we almost get one
correct digit.

\textbf{7. Concave partitions}

We have the generating function ( see [8] ) 
\[
f(x)=\prod_{j=1}^\infty (1-x^{j(j+1)/2})^{-1}=\sum_{n=0}^\infty a_nx^n 
\]
Observe that $a_n$ is also the number of partitions of $n$ into triangular
numbers. Then 
\[
g(t)=\log (f(\exp (-t))=\sum_{j=1}^\infty \sum_{k=1}^\infty \frac 1k\exp (-%
\frac{kj(j+1)t}2) 
\]
and the Mellin transformation 
\[
\widetilde{g}(s)=\Gamma (s)2^s\sum_{j=1}^\infty \sum_{k=1}^\infty \frac
1{k^{1+s}(j(j+1))^s} 
\]
\[
=2^s\Gamma (s)\zeta (1+s)H(s) 
\]
where 
\[
H(s)=\sum_{j=1}^\infty \frac 1{(j(j+1))^s} 
\]
We have to expand $H(s)$ near its poles and the poles of $\Gamma (s)\zeta
(1+s).$ Now 
\[
H(s)=\zeta (2s)+\sum_{j=1}^\infty \frac 1{j^{2s}}\left\{ (1+\frac
1j)^{-s}-1\right\} 
\]
\[
=\zeta (2s)+\sum_{k=1}^\infty \binom{-s}k\zeta (2s+k) 
\]
We want the residues at the singular points of 
\[
G(s)=2^s\Gamma (s)\zeta (1+s)t^{-s}\left\{ \zeta (2s)+\sum_{k=1}^\infty 
\binom{-s}k\zeta (2s+k)\right\} 
\]
The result is

\[
Res(G(s),s=\frac 12)=\frac{\sqrt{2\pi }}2\zeta (\frac 32)t^{-1/2} 
\]
\[
Res(G(s),s=0)=\log (\frac t{4\pi })-\gamma +\sum_{j=2}^\infty \frac{%
(-1)^j\zeta (j)}j=\log (\frac t{4\pi }) 
\]
since 
\[
\sum_{j=2}^\infty \frac{(-1)^j\zeta (j)}j=\gamma \text{ }=\text{Euler's
constant} 
\]
\[
Res(G(s),s=\frac 12-k)=\frac{\sqrt{2\pi }\zeta (\frac 12-k)t^{k+1/2}}{%
16\cdot 8^k(k+1)!} 
\]
We have 
\[
Res(G(s),s=-2k)=0 
\]
which is by no means trivial. It follows from the following identity:

\textbf{Proposition.}(Joakim Petersson). 
\[
\sum_{j=0}^k\binom kj\zeta (j-2k)=\frac{(-1)^k}{2\cdot (2k+1)\cdot \binom{2k}%
k} 
\]
Hence we have 
\[
f(\exp (-t))=\frac t{4\pi }\exp \left\{ \frac{\sqrt{2\pi }}2\zeta (\frac
32)\cdot t^{-1/2}+\frac{\sqrt{2\pi }}{16}\sum_{k=0}^\infty \frac{\zeta
(\frac 12-k)t^{k+1/2}}{8^k(k+1)!}\right\} 
\]
Put 
\[
a=\frac{\sqrt{2\pi }}2\zeta (\frac 32) 
\]
The main part of the first approximation of $a_n$ comes from 
\[
\frac t{4\pi }\exp (at^{-1/2})=\frac 1{4\pi }\sum_{j=0}^\infty \frac{%
a^jt^{1-j/2}}{j!}->\frac 1{4\pi }\sum_{j=0}^\infty \frac{a^jn^{j/2-2}}{%
j!\Gamma (j/2-1)}=q_0 
\]
(here we have $1/\Gamma (0)=0$). To get a more precise estimate we expand 
\[
\exp \left\{ \frac{\sqrt{2\pi }}{16}\sum_{k=0}^\infty \frac{\zeta (\frac
12-k)t^{k+1/2}}{8^k(k+1)!}\right\} =1+\sum_{k=1}^\infty c_kt^{k/2} 
\]
and get 
\[
q_k=\frac{c_k}{4\pi }\sum_{j=0}^\infty \frac{a^jn^{j/2-2-k/2}}{j!\Gamma
(j/2-1-k/2)} 
\]
We take a numerical example: Let $n=2000.$ Then 
\[
a_{2000}=27\text{ }79195\text{ }55391\text{ }39291 
\]
and 
\[
q_0=28\text{ }36990\text{ }09214\text{ }40082.213 
\]
\[
q_1=-58382\text{ }21952\text{ }56521.094 
\]
\[
q_2=594\text{ }72145\text{ }84992.877 
\]
\[
q_3=-8\text{ }07564\text{ }32951.404 
\]
\[
q_4=10135\text{ }44447.136 
\]
Thus 
\[
\Phi _1\thickapprox q_0+q_1+q_2+q_3+q_4= 
\]
\[
27\text{ }79194\text{ }61978\text{ }60047.307 
\]
so we get almost 7 correct digits.

\textbf{Remark.}

In [8] they show that 
\[
\log (a_n)\thicksim An^{1/3} 
\]
as $n\rightarrow \infty .$ Using the saddle point method and our $f(\exp
(-t))$ we obtain 
\[
A=3(\frac a2)^{2/3}=\frac 32\pi ^{1/3}\zeta (\frac 32)^{2/3} 
\]
\textbf{.}

\textbf{References.}

1. G.Almkvist, Partitions with parts in a finite set and with parts outside
a finite set, Exp. Math. 11 (2002),

449-456.

2. G. Almkvist, Asymptotic formulas and generalized Dedekind sums, Exp.
Math. 7 (1998), 343-359

3. G.Hardy,S.Ramanujan, Asymptotic formulae in combinatory analysis, Proc.
London Math. Soc.17

(1918), 75-115.

4. J. Knopfmacher, Abstract analytic number theory, Dover, N.Y. 1990

5. L.W.Kolitsch, M-order generalized Frobenius partitions with M colors, J.
Number Theory 39 (1991),

279-284.

6. I. Macdonald, Symmetric functions and Hall polynomials, 2-nd edition,
Oxford Science Publ.,1995.

7. J.M.Nolan,C.D.Savage,H. Wilf, Basic partitions, Preprint, Univ. of
Pennsylvania, 1995.

8. J.Snellman, M.Paulsen, Enumeration of concave integer partitions,
CO/0309065
\[
\]

Math Dept

Univ of Lund

Box 118

22100 Lund, Sweden

\end{document}